\documentclass[11pt]{amsart}
\usepackage{amsfonts,amscd}
\usepackage{amssymb}
\usepackage{url}

\theoremstyle{plain}
\newtheorem{theorem}                {Theorem}      [section]

\newtheorem{corollary}    [theorem]  {Corollary}

\theoremstyle{definition}

\newtheorem{remark}       [theorem]  {Remark}

\def \r{{\mathbb R}}
\def \s{{\mathbb S}}
\def \z{{\mathbb Z}}
\def \c{{\mathbb C}}

\def \hor{\mathcal H}
\def \b{\mathcal B}

\numberwithin{equation}{section}

\begin{document}

\title[A generalisation of the Hopf Construction]{A generalisation of the Hopf Construction and harmonic morphisms into $\s^2$}

\author{S.~Montaldo}
\address{Universit\`a degli Studi di Cagliari\\
Dipartimento di Matematica e Informatica\\
Via Ospedale 72\\
09124 Cagliari, Italia}
\email{montaldo@unica.it}

\author{A.~Ratto}
\address{Universit\`a degli Studi di Cagliari\\
Dipartimento di Matematica e Informatica\\
Viale Merello 93\\
09123 Cagliari, Italia}
\email{rattoa@unica.it}

\begin{abstract}
In this paper we construct a new family of harmonic morphisms $\varphi:V^5\to\s^2$, where $V^5$ is
a $5$-dimensional open manifold contained in an ellipsoidal hypersurface of
$\c^4=\r^8$. These harmonic morphisms admit a continuous extension to the
completion ${V^{\ast}}^5$, which turns out to be an explicit real algebraic variety.
We work in the context of a generalization of the Hopf construction and equivariant theory.
\end{abstract}

\date{}



\thanks{}

\maketitle

\section{Introduction}
{\it Harmonic} maps  are critical points of the {\em energy} functional
$$
E(\varphi)=\frac{1}{2}\int_{M}\,|d\varphi|^2\,dv_g,
$$
where $\varphi:(M,g)\to(N,h)$ is a smooth map between two Riemannian
manifolds $M$ and $N$.
This is a very wide area of research, involving a rich interplay of geometry,
analysis and topology. We refer to \cite{JEELLE1,JEELLE2} for notation and
background on harmonic maps and to \cite{BMBib} fore a more recent
bibliography.

A geometrically significant sub-family of harmonic maps is represented by the so-called harmonic morphisms.
An exhaustive reference for this topic is the book of P.~Baird and
J.C.~Wood \cite{PBJCW}, where characterizing properties and existence of
harmonic morphisms are presented in connection with central themes such as harmonic functions and potential
theory, together with conformal mappings in the plane and holomorphic maps into a Riemann surface.

For an operational point of view the simplest way to characterize harmonic morphisms is to say that
they are just harmonic maps with the additional property that the differential
$d\varphi$ is a horizontally weakly conformal map: this means that, at any point $x\in M$,
either $d\varphi_x$ vanishes or
\begin{equation}\label{eq:hor-conf}
d\varphi_x:(T_x M)_{\hor}\to T_{\varphi(x)}N
\end{equation}
is surjective and conformal. More precisely, in \eqref{eq:hor-conf},
$(T_x M)_{\hor}$ denotes the horizontal space
$(\ker(d\varphi_x))^{\perp}$ and it is required that there exists a
number $\Lambda(x)>0$ such that
$$
\Lambda(x) g(X,Y)= h(d\varphi_x(X),d\varphi_x(Y)),\quad \forall
X,Y\in (T_x M)_{\hor}
$$
The function $\lambda(x)=\sqrt{\Lambda(x)}$ is called the {\it
dilation} (of $\varphi$ at $x$). In particular, if $\varphi$ is a
non constant harmonic morphism, then $m\geq n$, where
 $m=\dim(M)$ and $n=\dim(N)$. Moreover, the set of singular points
(i.e., those points where $d\varphi_x=0$) is a closed {\it polar set}, i.e., it has
zero capacity (see \cite{PBJCW} for details).

In this paper we work in the context of equivariant theory: roughly speaking,
that means that we restrict our attention to a class of mappings having enough
symmetries  to guarantee that harmonicity reduces to the study of a second order {\it ordinary}
differential equation.
We refer to \cite{EELRAT} for notation, background and examples.
More specifically, we shall propose a generalisation of the
Hopf construction which gives rise to a new family
of harmonic morphisms from a $5$-dimensional manifold with
singularities onto the Euclidean $2$-sphere $\s^2$.

\section{Statement of the main results and related comments}

In order to illustrate our framework, we first consider the $3$-dimensional
ellipsoid
\begin{equation}\label{eq:ellipsoid}
{Q^{\ast}}^3(a,b)=\left\{[x,y]\in\c\times\c\;:\;\frac{|x|^2}{a^2}+\frac{|y|^2}{b^2}=1\right\},\quad (a,b>0)
\end{equation}
In \cite[chapter X]{EELRAT} the authors study in detail a family
of maps
\begin{equation}\label{eq:fkl}
\varphi_{k,\ell}:{Q^{\ast}}^3(a,b)\to \s^2\subset\r^2\times\r,\quad k,\ell\in\z
\end{equation}
of the following form (Hopf's construction)
\begin{equation}\label{eq:hopf}
[a \sin s\;e^{i\theta_1},b \cos s\;e^{i\theta_2}]\mapsto
[\sin\alpha(s) \;e^{i(k\theta_1+\ell\theta_2)}, \cos\alpha(s)]
\end{equation}
where the function $\alpha:(0,\pi/2)\to(0,\pi)$
satisfies the boundary conditions
\begin{equation}\label{eq:boundary}
(i)\:\lim_{s\to0^+}\alpha(s)=0,\quad(ii)\:\lim_{s\to\tfrac{\pi}{2}^-}\alpha(s)=\pi.
\end{equation}
In particular, it turns out that the map $\varphi_{k,\ell}$ is harmonic if and only if $\alpha$
 is a solution of the ODE
$$
\alpha''+(\cot s-\tan s)\alpha'-\frac{h'(s)}{h(s)}\alpha'-h^2(s)\left(\frac{k^2}{a^2 \sin^2 s}+\frac{\ell^2}{b^2 \cos^2 s}\right)\sin \alpha \cos \alpha=0
$$
where $h^2(s)=a^2 \cos^2 s+b^2 \sin^2 s$. Note, for future
comparison, that the left member of \eqref{eq:hopf} represents a
parametrisation of ${Q^{\ast}}^3(a,b)$ ($0\leq s \leq {\pi}/{2},\;
0\leq\theta_1,\theta_2<2\pi$). The boundary conditions
\eqref{eq:boundary} ensure that the solutions have a regular
extension through the loci $s=0$ and $s=\pi/2$, where the
coefficients of the ODE become singular. The map $\varphi_{k,\ell}$
has topological significance because it has Hopf invariant $k\ell$,
i.e., it represents $k\ell\in\z=\pi_3(\s^2)$. In \cite{EELRAT} it is
proved that, under the additional restriction
\begin{equation}\label{eq:aoverb}
\frac{a}{b}=\left|\frac{\ell}{k}\right|,
\end{equation}
there exist harmonic morphisms of the type \eqref{eq:hopf}.

We are now in the right position to proceed to our generalisation of the previous examples.

We shall consider ${Q^{\ast}}^5={Q^{\ast}}^5(a_1,a_2,a_3,a_4)\subset\c^4,\;a_1,a_2,a_3,a_4>0$,
where ${Q^{\ast}}^5$ is the subset of $\c^4$ parametrised by
\begin{equation}\label{eq:paramtr-q5}
[a_1 \sin s\;e^{i\theta_1},a_2 \sin(s+\tfrac{\pi}{4})\;e^{i\theta_2},a_3 \cos s\;e^{i\theta_3},a_4 \cos(s+\tfrac{\pi}{4})\;e^{i\theta_4}]
\end{equation}
with $0\leq\theta_1,\theta_2,\theta_3,\theta_4<2\pi$, $0\leq s\leq \pi/4$. Away from the loci $s=0$ and $s=\pi/4$, \eqref{eq:paramtr-q5} represents a Riemannian manifold
\begin{equation}\label{eq:q5-riemannian}
(Q^5,g),
\end{equation}
where $Q^5=\s^1\times\s^1\times\s^1\times\s^1\times(0,\pi/4)$, and the Riemannian metric
$g$ induced by the Euclidean structure of $\c^4(=\r^8)$ is given by
\begin{align}\label{eq:q5-metric}
g=&[a_1^2 \sin^2s]\; d\theta_1^2+[a_2^2 \sin^2(s+\tfrac{\pi}{4})]\;
d\theta_2^2+[a_3^2 \cos^2s]\; d\theta_3^2 \nonumber\\
&+[a_4^2 \cos^2(s+\tfrac{\pi}{4})]\; d\theta_4^2+h^2(s) ds^2,
\end{align}
with
\begin{equation}\label{eq:q5-h}
h(s)=\sqrt{a_1^2 \cos^2s+a_2^2 \cos^2(s+\tfrac{\pi}{4})
+a_3^2 \sin^2s+a_4^2 \sin^2(s+\tfrac{\pi}{4})}
\end{equation}

\begin{remark}
The locus $s=0$ is not a topological singularity, since across it ${Q^{\ast}}^5$ is
locally homeomorphic to $U^2\times\s^1\times\s^1\times\s^1$, where $U^2$ is an open subset of $\r^2$.
However, performing the change of local coordinates
$$
x_1=a_1 \sin s \cos \theta_1,\;x_2=a_1 \sin s \sin \theta_1,\;
\theta_2=\theta_2,\;\theta_3=\theta_3,\;\theta_4=\theta_4
$$
it is not difficult to check that the coefficients of the metric
tensor $g$ are continuous, but not differentiable across $s=0$. The
same happens across $s=\pi/4$.
\end{remark}

Bearing in mind the family of maps $\varphi_{k,\ell}$ of
\eqref{eq:fkl}, we now define the class of equivarinat mappings we are interested in ($k_i\in\z,\,i=1,\ldots,4$):
$$
\varphi=\varphi_{k_1,k_2,k_3,k_4}:Q^5(a_1,a_2,a_3,a_4)\to\s^2\subset\r^2\times\r
$$
where, to simplify notation, we write $w$ in place of the generic point \eqref{eq:paramtr-q5} of $Q^5$ and require

\begin{equation}\label{eq:def-varphi}
\varphi(w)=[\sin\alpha(s)\; e^{i(k_1\theta_1+k_2\theta_2+k_3\theta_3+k_4\theta_4)},\cos\alpha(s)].
\end{equation}

Here the function $\alpha:(0,\pi/4)\to(0,\pi)$ must satisfy the boundary conditions (compare with \eqref{eq:boundary}):

\begin{equation}\label{eq:boundary-2}
(i)\:\lim_{s\to0^+}\alpha(s)=0,\quad(ii)\:\lim_{s\to\pi/4^-}\alpha(s)=\pi.
\end{equation}

Such boundary conditions ensure that the map $\varphi$ as in \eqref{eq:def-varphi} extends {\it continuously}
across the loci $s=0$ and $s=\pi/4$ in ${Q^{\ast}}^5$. Harmonicity of such $\varphi$ depends on an ODE for $\alpha$:
we shall prove that this ordinary differential equation, under suitable restrictions of the type \eqref{eq:aoverb},
admits a prime integral which turns out to be essentially equivalent to the horizontal conformality of the map.
A qualitative study of this prime integral will lead us to the existence of strictly increasing solutions satisfying
\eqref{eq:boundary-2}. More precisely, our result is

\begin{theorem}\label{teo-main}
Assume that $a_i=|k_i|,\, i=1,\ldots,4$. Then there exist harmonic morphisms $\varphi:Q^5\to\s^2$ of the
type \eqref{eq:def-varphi} which verify the boundary conditions \eqref{eq:boundary-2}, so admitting a continuous
extension to the whole ${Q^{\ast}}^5$.
\end{theorem}

In view of the previous theorem, it is natural to study the manifold
${Q^{\ast}}^5$ a bit more accurately. First of all, we observe
(routine verification using \eqref{eq:paramtr-q5}) that all the
points of
  ${Q^{\ast}}^5$ satisfy the following set of polynomial equations
\begin{equation}\label{eq:q5equations}
\left\{\begin{array}{lr}
\displaystyle{\frac{|z_1|^2}{a^2_1}+\frac{|z_3|^2}{a^2_3}=1}&(i)\\
\displaystyle{\frac{|z_2|^2}{a^2_2}+\frac{|z_4|^2}{a^2_4}=1}&(ii)\\
\displaystyle{\left[ \frac{|z_3|^2}{a^2_3}-\frac{|z_1|^2}{a^2_1}\right]^2+\left[ \frac{|z_4|^2}{a^2_4}-\frac{|z_2|^2}{a^2_2}\right]^2=1}&(iii)\\
\end{array}
\right.
\end{equation}
where $[z_1,z_2,z_3,z_4]$ denote complex coordinates on $\c^4$. Next, let us call
${V^{\ast}}^5={V^{\ast}}^5(a_1,a_2,a_3,a_4)$ the real algebraic variety defined by the set of equations \eqref{eq:q5equations}.

We have just remarked that ${Q^{\ast}}^5\subset{V^{\ast}}^5$. On the other hand, a combined inspection of \eqref{eq:paramtr-q5}
 and \eqref{eq:q5equations} shows that a parametrisation of the whole ${V^{\ast}}^5$ requires, in  \eqref{eq:paramtr-q5},
 to let $s$ vary in the interval $[0,\pi]$. Now the singular loci are $s=0,\,s=\pi/4,\,s=\pi/2$ and $s=3\pi/4$
 ($s=\pi$ coincide with $s=0$). Away from these loci, we have a Riemannian manifold $(V^5,g)$ as follows (compare with \eqref{eq:q5-riemannian}):
\begin{equation}\label{eq:v5-riemannian}
V^5=\s^1\times\s^1\times\s^1\times\s^1\times I
\end{equation}
where
$$
I=(0,\tfrac{\pi}{4})\cup (\tfrac{\pi}{4},\tfrac{\pi}{2})\cup (\tfrac{\pi}{2},\tfrac{3\pi}{4})\cup (\tfrac{3\pi}{4},\pi)
$$
 and the metric $g$ is as in
\eqref{eq:q5-metric}, with $s\in I$.

Next, it becomes natural to study maps
\begin{equation}\label{eq:varphi-v5}
\varphi=\varphi_{k_1,k_2,k_3,k_4}:V^5(a_1,a_2,a_3,a_4)\to\s^2\subset\r^2\times\r
\end{equation}
defined again by \eqref{eq:def-varphi}, with the only difference that this time the function $\alpha:I\to (0,4\pi)$
satisfies the following set of boundary conditions

\begin{equation}\label{eq:boundary3}
\begin{array}{lll}
(i)&\displaystyle{\lim_{s\to 0^+}\alpha(s)=0},&\displaystyle{\lim_{s\to \tfrac{\pi}{4}^-}\alpha(s)=\pi}\\
(ii)&\displaystyle{\lim_{s\to \tfrac{\pi}{4}^+}\alpha(s)=\pi},&\displaystyle{\lim_{s\to \tfrac{\pi}{2}^-}\alpha(s)=2\pi}\\
(iii)&\displaystyle{\lim_{s\to \tfrac{\pi}{2}^+}\alpha(s)=2\pi},&\displaystyle{\lim_{s\to \tfrac{3\pi}{4}^-}\alpha(s)=3\pi}\\
(iv)&\displaystyle{\lim_{s\to \tfrac{3\pi}{4}^+}\alpha(s)=3\pi},&\displaystyle{\lim_{s\to \pi^-}\alpha(s)=4\pi}.
\end{array}
\end{equation}
\begin{remark}
The boundary conditions \eqref{eq:boundary3} ensure that the maps $\varphi$ in
\eqref{eq:varphi-v5} extend continuously to the whole ${V^{\ast}}^5$. In particular, such extension of the
$\varphi$'s cover $4$ times the $2$-sphere $\s^2$: the singular loci $s=0$ and $s=\pi/2$ are sent to the North pole, while
$s=\pi/4$ and $s=3\pi/4$ go into the South pole.
\end{remark}

An adaptation of the proof of Theorem~\ref{teo-main} enables us to obtain the following

\begin{corollary}\label{cor}
Assume that $a_i=|k_i|$, $i=1,\ldots,4$. Then there exist harmonic morphisms $\varphi:V^5\to\s^2$ of the type
\eqref{eq:varphi-v5} which satisfy the boundary conditions
\eqref{eq:boundary3}, so admitting a continuous extension to the whole ${V^{\ast}}^5$.
\end{corollary}

\section{Proof of the results and final remarks}

{\em Proof of Theorem~\ref{teo-main}.} We have to show that, under the assumption $a_i=|k_i|,\, i=1,\ldots,4$,
there exist a function $\alpha:(0,\pi/4)\to(0,\pi)$ which satisfies the boundary conditions  \eqref{eq:boundary-2}
 and such that the associated $\varphi$
of \eqref{eq:def-varphi} is both harmonic and horizontally
conformal. To this we need to express, explicitly, conditions which
are equivalent to harmonicity and horizontal conformality.

By applying, for instance, the Reduction Theorem~(4.13) of \cite{EELRAT}, it is easy
to verify that the harmonicity equation for maps $\varphi$ as in \eqref{eq:def-varphi} is given by
\begin{equation}\label{eq:ode-alpha}
\alpha''(s)+D(s)\alpha'(s)-G(s,\alpha,\alpha')=0,
\end{equation}
where the functions $D$ and $G$ can be calculated as in \cite[p.
153]{EELRAT}. One finds
\begin{equation}\label{eq:ds}
D(s)=\left[ \frac{\cos s}{\sin s}+\frac{\cos (s+\tfrac{\pi}{4})}{\sin (s+\tfrac{\pi}{4})}-\frac{\sin s}{\cos s}
-\frac{\sin (s+\tfrac{\pi}{4})}{\cos(s+\tfrac{\pi}{4})}\right]-\frac{h'(s)}{h(s)}
\end{equation}
where $h(s)$ is the function defined in \eqref{eq:q5-h}. And also

\begin{align} \label{eq:gs}
G(s,\alpha,\alpha')=&h^2(s)\left[ \frac{k_1^2}{a_1^2 \sin^2 s}+\frac{k_2^2}{a_2^2 \sin^2(s+\tfrac{\pi}{4})}+\right.\nonumber\\
&\left.+\frac{k_3^2}{a_3^2 \cos^2 s}+\frac{k_4^2}{a_4^2
\cos^2(s+\tfrac{\pi}{4})}\right]\sin\alpha(s) \cos\alpha(s).
\end{align}
Now we simplify \eqref{eq:ds} and \eqref{eq:gs}. Indeed, using three times the trigonometrical identity
\begin{equation}\label{eq:trig}
\cot x -\tan x = 2\cot 2x
\end{equation}
we find that \eqref{eq:ds} ca be rewritten as
\begin{equation}\label{eq:ds-simplfy}
D(s)=4 \cot(4s)-\frac{h'(s)}{h(s)}.
\end{equation}
Next, we use the hypothesis $a_i=|k_i|,\, i=1,\ldots,4$, in \eqref{eq:gs}. A routine
computation then leads us to
\begin{equation}\label{eq:gs-simplfy}
G(s,\alpha,\alpha')=h^2(s) \frac{16}{\sin^2(4s)} \sin\alpha(s) \cos\alpha(s).
\end{equation}
By way of summary, replacing \eqref{eq:ds-simplfy} and \eqref{eq:gs-simplfy} in
\eqref{eq:ode-alpha}, we have proved that the harmonicity equation for the maps $\varphi$ as in  \eqref{eq:def-varphi} is given by
\begin{equation}\label{eq:ode-alpha-simplify}
\alpha''(s)+\left[4 \cot(4s)-\frac{h'(s)}{h(s)}\right]\alpha'(s)-\frac{16 h^2(s)}{\sin^2(4s)} \sin\alpha(s) \cos\alpha(s)=0.
\end{equation}

The next step is to work out the condition for horizontal conformality. An orthonormal base of the tangent space $T_wQ^5$ is
$$
\b=\{e_1,e_2,e_3,e_4,e_5\}
$$
where
\begin{equation}\label{eq:base}
\begin{array}{lll}
\displaystyle{e_1=\frac{1}{a_1\sin
s}\frac{\partial}{\partial\theta_1}}&
\displaystyle{e_2=\frac{1}{a_2\sin(s+\tfrac{\pi}{4})}\frac{\partial}{\partial\theta_2}}&
\displaystyle{e_3=\frac{1}{a_3\cos s}\frac{\partial}{\partial\theta_3}}\\
&\\
\displaystyle{e_4=\frac{1}{a_4\cos(s+\tfrac{\pi}{4})}\frac{\partial}{\partial\theta_4}}&
\displaystyle{e_5=\frac{1}{h(s)}\frac{\partial}{\partial s}}.&
\end{array}
\end{equation}
We write the Riemannian metric on the range $\s^2$ as
\begin{equation}\label{eq:metrics2}
h=\sin^2t\; d\gamma^2+dt^2
\end{equation}
since in  \eqref{eq:def-varphi} we are parametrising $\s^2$, as a
subset of $\r^2\times\r$, by means of
\begin{equation}\label{eq:parametrisations2}
[\sin t\; e^{i\gamma}, \cos t]\quad t\in[0,\pi],\;\gamma\in[0,2\pi).
\end{equation}
Now we can compute the differential for maps $\varphi$ as in \eqref{eq:def-varphi}.
We find
\begin{equation}\label{eq:diffvarphi}
\begin{array}{l}
\displaystyle{d\varphi(e_1)=\frac{k_1}{a_1\sin s}
\frac{\partial}{\partial\gamma}=\frac{1}{\sin s}\frac{\partial}{\partial\gamma}}\\
\\
\displaystyle{d\varphi(e_2)=\frac{k_2}{a_2\sin(s+\tfrac{\pi}{4})}
\frac{\partial}{\partial\gamma}=\frac{1}{\sin(s+\tfrac{\pi}{4})}\frac{\partial}{\partial\gamma}}\\
\\
\displaystyle{d\varphi(e_3)=\frac{k_3}{a_3\cos s}
\frac{\partial}{\partial\gamma}=\frac{1}{\cos s}\frac{\partial}{\partial\gamma}}\\
\\
\displaystyle{d\varphi(e_4)=\frac{k_4}{a_4\cos(s+\tfrac{\pi}{4})}
\frac{\partial}{\partial\gamma}=\frac{1}{\cos(s+\tfrac{\pi}{4})}\frac{\partial}{\partial\gamma}}\\
\\
\displaystyle{d\varphi(e_5)=\frac{\alpha'(s)}{h(s)}
\frac{\partial}{\partial t}},\\
\end{array}
\end{equation}
where we have used $a_i=k_i,\, i=1,\ldots,4$ ( the case $a_i=-k_i$ can be handled similarly and it is left to
 the interested reader). Now, let $v=\sum_{i=1}^5 v_i e_i$ be
a generic vector in $T_w Q^5$. By using \eqref{eq:diffvarphi} and the linearity of $d\varphi$ we find that
\begin{equation}\label{eq:vker}
v\in \ker(d\varphi) \Leftrightarrow
\begin{cases}
v_5=0 \\
\displaystyle{
\frac{v_1}{\sin s}+\frac{v_2}{\sin(s+\tfrac{\pi}{4})}+\frac{v_3}{\cos s}+\frac{v_4}{\cos(s+\tfrac{\pi}{4})}=0}
\end{cases}.
\end{equation}
Therefore, we deduce that a base for the horizontal part is
\begin{equation}\label{eq:hor-part}
\{y,e_5\}\quad\text{with}\quad y=\frac{e_1}{\sin s}+\frac{e_2}{\sin(s+\tfrac{\pi}{4})}+\frac{e_3}{\cos s}+\frac{e_4}{\cos(s+\tfrac{\pi}{4})}.
\end{equation}

Since

$$
\|y\|^2=\frac{1}{\sin^2
s}+\frac{1}{\sin^2(s+\tfrac{\pi}{4})}+\frac{1}{\cos^2
s}+\frac{1}{\cos^2(s+\tfrac{\pi}{4})}=\frac{16}{\sin^2(4s)}
$$
we conclude that an orthonormal base for $(T_wQ^5)_{\hor}$ is
$$
\b_{\hor}=\{y^{\ast},e_5\}
$$
where
$$
y^{\ast}=\frac{y}{\|y\|}=\frac{\sin(4s)}{4} \;y.
$$
Because $d\varphi$ preserves orthogonality between $y^{\ast}$ and $e_5$, horizontal conformality reduces to the condition
\begin{equation}\label{eq:confye5}
\|d\varphi(y^{\ast})\|^2=\|d\varphi(e_5)\|^2.
\end{equation}
Now, using \eqref{eq:base} and \eqref{eq:diffvarphi}, we easily find
\begin{eqnarray}\label{eq:dvarphiystar}
d\varphi(y^{\ast})&=&\left[\frac{1}{\sin^2 s}+\frac{1}{\sin^2(s+\tfrac{\pi}{4})}+\frac{1}{\cos^2 s}+\frac{1}{\cos^2(s+\tfrac{\pi}{4})} \right]\frac{\sin(4s)}{4}\frac{\partial}{\partial\gamma}\nonumber\\
&=&\frac{4}{\sin(4s)}\frac{\partial}{\partial\gamma}.
\end{eqnarray}
Consequently, we have
\begin{equation}\label{eq:normdvarphiystar}
\|d\varphi(y^{\ast})\|^2=\frac{16}{\sin^2(4s)}\left\|\frac{\partial}{\partial\gamma}\right\|^2=
\frac{16}{\sin^2(4s)}\; \sin^2\alpha(s).
\end{equation}
Moreover, using \eqref{eq:diffvarphi}
\begin{equation}\label{eq:normdvarphie5}
\|d\varphi(e_5)\|^2=\frac{{\alpha'}^2(s)}{h^2(s)}\left\|\frac{\partial}{\partial t}\right\|^2=
\frac{{\alpha'}^2(s)}{h^2(s)}.
\end{equation}
By way of summary, putting together \eqref{eq:confye5},
\eqref{eq:normdvarphiystar} and \eqref{eq:normdvarphie5} we obtain
that horizontal conformality is equivalent to
\begin{equation}\label{eq:conformalityfinall}
\frac{16 \sin^2\alpha(s)}{\sin^2(4s)}=
\frac{{\alpha'}^2(s)}{h^2(s)}.
\end{equation}
We are now in the right position to end the proof of Theorem~\ref{teo-main}. First, taking square roots and
re-arranging the terms, we rewrite \eqref{eq:conformalityfinall} in a more convenient form
\begin{equation}\label{eq:conformalityfinalleasier}
\frac{\alpha'(s)}{\sin\alpha(s)}=
\frac{4h(s)}{\sin(4s)}.
\end{equation}
Now, applying $\frac{d}{d s}$ to both members of \eqref{eq:conformalityfinalleasier}
and computing, one finds that \eqref{eq:conformalityfinalleasier} implies
\eqref{eq:ode-alpha-simplify}, or, to say it in words, the condition  \eqref{eq:conformalityfinalleasier} of
horizontal conformality is a prime integral for the harmonicity equation \eqref{eq:ode-alpha-simplify}.
Thus it is enough to show that the first order
ODE \eqref{eq:conformalityfinalleasier} admits a solution $\alpha:(0,\tfrac{\pi}{4})\to(0,\pi)$ which
satisfies the prescribed boundary conditions
\eqref{eq:boundary-2}.
Indeed, \eqref{eq:conformalityfinalleasier} can be integrated by separation of variables and that yields
 a family of solutions
\begin{equation}\label{eq:solutions}
\alpha(s)=2 \tan^{-1}\left[ c \exp\left( \int_{\frac{\pi}{8}}^s 4
\frac{h(u)}{\sin(4u)}du\right)\right],\quad s\in(0,\tfrac{\pi}{4}),
\;c>0.
\end{equation}
Since the positive function $h(u)$ is bounded and bounded away from zero, mere
inspection of the solutions \eqref{eq:solutions} enables us to conclude that the boundary conditions
\eqref{eq:boundary-2} are fulfilled, so ending the proof of the theorem.

\begin{remark}
The explicit meaning of the integration constant $c$ in  \eqref{eq:solutions} is
$$
c=\left|\tan\left(\frac{\alpha(\tfrac{\pi}{8})}{2}\right)\right|.
$$
In particular, the most symmetric solution, i.e. the one with $\alpha(\tfrac{\pi}{8})=\tfrac{\pi}{2}$,
 corresponds to $c=1$.
\end{remark}

{\em Proof of Corollary~\ref{cor}.} Proceeding precisely as in the proof of Theorem~\ref{teo-main},
 we have again that the horizontal conformality
 \eqref{eq:conformalityfinalleasier} is a prime integral for the
 harmonicity equation \eqref{eq:ode-alpha-simplify}. Now one needs
 to show the existence of a solution $\alpha:I\to(0,4\pi)$ which
 satisfies the boundary conditions \eqref{eq:boundary3}. The
 construction of such $\alpha$ is divided into $4$ Steps: in Step $1$,
 one obtains $\alpha:(0,{\pi}/{4})\to(0,\pi)$ which satisfies
\eqref{eq:boundary3} $(i)$. Of course, this was done during the
proof of Theorem~\ref{teo-main}. Step $2$ is the construction of a
solution $\alpha:({\pi}/{4},{\pi}/{2})\to(\pi,2\pi)$
satisfying \eqref{eq:boundary3} $(ii)$. That is achieved again by
explicit integration of \eqref{eq:conformalityfinalleasier} which
gives ($c>0$)
\begin{equation}\label{eq:solutions-2}
\alpha(s)=2\pi+2 \tan^{-1}\left[ -c \exp\left(
\int_{\frac{3\pi}{8}}^s 4
\frac{h(u)}{\sin(4u)}du\right)\right],\quad
s\in(\tfrac{\pi}{4},\tfrac{\pi}{2}).
\end{equation}
Once more, direct inspection of \eqref{eq:solutions-2} confirms that
the boundary conditions \eqref{eq:boundary3} $(ii)$ are verified. In
a similar fashion, one constructs solutions
$\alpha:({\pi}/{2},{3\pi}/{4})\to(2\pi,3\pi)$ satisfying
 \eqref{eq:boundary3} $(iii)$, and finally $\alpha:({3\pi}/{4},\pi)\to(3\pi,4\pi)$ which verify
 \eqref{eq:boundary3} $(iv)$, so completing the proof of the
 corollary.

 \begin{remark}
 It follows easily from \eqref{eq:q5equations} $(i)$ and $(ii)$ that
 ${V^{\ast}}^5$ is contained into the ellipsoidal hypersurface of
 equation
\begin{equation}\label{eq:ellipsoid-z}
\frac{|z_1|^2}{a_1^2}+\frac{|z_2|^2}{a_2^2}+\frac{|z_3|^2}{a_3^2}+\frac{|z_4|^2}{a_4^2}=2.
\end{equation}
Under the additional restriction
\begin{equation}\label{eq:restrictions}
a_1=a_3\quad\text{and}\quad a_2=a_4
\end{equation}
the function $h(s)$ of \eqref{eq:q5-h} is constant, i.e.,
\begin{equation}\label{eq:h-constant}
h(s)\equiv\sqrt{a_1^2+a_2^2}=A.
\end{equation}
In this case the solutions \eqref{eq:solutions} of
Theorem~\ref{teo-main} (and, similarly, the solutions of
Corollary~\ref{cor}) can be made more explicit. For instance,
performing an integration leads us to express \eqref{eq:solutions},
under the hypothesis \eqref{eq:restrictions}, as
\begin{equation}\label{eq:solutions-3}
\alpha(s)=2\tan^{-1}\left[ c\, (\tan 2s)^A\right],\quad
s\in(0,\tfrac{\pi}{4}).
\end{equation}
\end{remark}
\begin{remark}
The integration constant $c>0$ in \eqref{eq:solutions} reflects the
fact that our construction provides an infinite family of solutions.
A similar situation had been observed in \cite[p. 186]{EELRAT}, and
\cite[Chapter 13]{PBJCW}, where prime integrals occur and produce
variations through equivariant harmonic morphisms. It should also
be observed that all the solutions in Theorem~\ref{teo-main} and
Corollary~\ref{cor} have $\alpha'(s)>0$: the simplest way to verify
this claim is direct inspection of
\eqref{eq:conformalityfinalleasier}. In particular, only the singular loci are sent into the poles of $\s^2$.
\end{remark}
\begin{remark}
For the purpose of comparison, we also point out that some other
examples of harmonic morphisms defined on open subsets of compact
manifolds can be found in \cite[p. 410]{PBJCW}.
\end{remark}
\begin{remark}
Since the range of our harmonic morphisms is $2$-dimensional, their
fibres provide a foliation of $V^5$ made of $3$-dimensional minimal
submanifolds diffeomorphic to $\s^1\times\s^1\times\s^1$.
\end{remark}

\end{document}